\input amstex
\documentstyle{amsppt}
\document
\magnification=1200 
\NoBlackBoxes
\vsize16cm
\hoffset=0.7in
\nologo

%\hfill{\it file Amywork/balzan.tex, version 21.05.2008}

%\medskip

\centerline{\bf TRUTH AS VALUE AND DUTY:}

\smallskip 

\centerline{\bf LESSONS OF MATHEMATICS\footnotemark1}
\footnotetext{Talk at the International  Symposium of
the Balzan Foundation  ``Truth in the Humanities, Science and Religion'',
Lugano, May 16--17,  2008.}

\bigskip
\centerline{\bf Yu.~I.~Manin}

\medskip

\centerline{\it Max--Planck--Institut f\"ur Mathematik, Bonn, Germany,}

\centerline{\it and Northwestern University, Evanston, USA}

\bigskip

\centerline{\bf 1. Introduction}

\medskip

Imagine that you open your morning newspaper and read the following report:

\medskip

 {\bf Brownsville, AR.} {\it A local object partially immersed in a liquid
was buoyed upward Tuesday by a force equal to the weight
of the liquid displaced by that object, witnesses
at the scene reported. As of press time, the object is still maintaining
positive buoyancy.}

\medskip

In fact, I did read this  report in {\it the ONION}; I have only  
abridged it to add a F\'en\'eonian touch.

\smallskip

If this meeting had been dedicated to the nature of the comical,
one could produce an interesting analysis of the clever
silliness of this parody. But as we are preoccupied
with truth, I will use it in order to illustrate
the differences between the attitudes to truth among practitioners
of social sciences and law
(as exemplified by [BoHa]) on the one hand, and that
of, say, physicists, on the other. 

\smallskip

To put it crudely, in social sciences information comes from witnesses;
but in what sense was Archimedes' role in his discovery  that of a witness,
and are the experimental observations generating/supporting a physical theory 
on an equal footing with the observations of witnesses to a crime scene,
or respondents to a poll?

\smallskip

Now, imagine another report, that could have been posted on the web--site of the
Department of Physics of Cambridge University:

\medskip

{\it The  Cavendish Laboratory News \& Features bulletin announced yesterday
that a Cavendish student has won Science, Engineering and Technology award.
He  managed to measure the constant $\pi$
with unprecedented precision: $\pi = 3,1415925 ...$ with 
an error $\pm 2$ at  the last digit.}

\medskip

I must confess right away that I did not
read  but simply fabricated this spoof
in order to stress the further differences between the attitudes towards  truth, now held by
physicists  and by mathematicians respectively. 

\smallskip

On the one hand, formally such an announcement would make perfect sense:
the {\it mathematical} constant $\pi$ {\it can} be measured with
some precision, in the same way that any {\it physical}
constant such as the speed of light $c$, or the mass of the electron can be measured.
The maximum achievable precision, at least of a ``naive''
direct measurement of $\pi$, is determined by the degree to which 
we can approximate ideal Euclidean rigid bodies by real
physical ones. The limits to this approximation are set
by the atomic structure of matter, and in the final
analysis, by quantum effects. 

\smallskip

On the other hand, in order to get in principle
as many digits of $\pi$ as one wishes,  
measurements are not required at all. Instead, one can use one of the many existing
formulas/algorithms/software codes and do it
on a sheet of paper, a pocket calculator, or a supercomputer.
This time the limits of precision are determined by the
physical limitations of our calculator: the size
of the sheet of paper, memory of computer, construction
of the output device, available time ... 

\smallskip

What I want to stress now is that $\pi$ imagined as
an infinite sequence of its digits,  is not amenable
to a ``finite'' calculation: even the number of digits of $\pi$
equal to the number of atoms in the observable Universe,
would not exhaust $\pi$. Nevertheless, mathematicians speak about $\pi$
and work with $\pi$  as if it were a 
completely well defined entity, graspable in its entirety
not only by one exceptional super--Mind, but by the minds
of all trained researchers, never doubting that
when they speak of $\pi$, they speak about one and the same
ideal object, as rigid as if it really exists  in some Platonic world.

\smallskip

In fact, one facet of this rigidity can be expressed
by a few theorems implying that whatever exact formula, algorithm,
or software code we might use to calculate $\pi$ and whatever
precision we choose, we will always get the same result.
If we do not, either our formula was wrong, or the calculator
made a mistake/there was a bug in the code/output device
could not cope with the quantity of information ...

\smallskip

Contemplating this example, we may grasp 
the meaning of the succinct description
of mathematics by Davis and Hersh ([DaHe]): {\it ``the study of mental
objects with reproducible properties''.}  
 
\smallskip

However, I want to use this example in order to stress that most of 
the deep mathematical truths are about infinity and 
infinitary mental constructs rather than experimentally verifiable
finitary -- and finite -- operations, that can be modeled using actual objects of the physical world.

\newpage

\centerline{\bf 2}

\medskip

\hfill{\it ... mais je ne le crois pas!}

\smallskip

\hfill{\it G.~Cantor to R.~Dedekind, June 29, 1877}

\medskip

Before Georg Cantor, infinity appeared in mathematical
theorems mostly implicitly, through the quantifier ``all'' (which also
could be only implicit as in most Euclid's theorems). 

\smallskip

Cantor proved the first  theorem ever in which infinities themselves were
objects of consideration. Slightly modernizing his arguments,
we can say that he invented two or three mental constructions
allowing us to compare sizes (technically, cardinalities) 
of infinite (in fact, finite as well) sets:

\smallskip

a) Two sets $X$, $Y$ have equal cardinalities, symbolically $|X|=|Y|$, 
if their elements $x\in X$, $y\in Y$ can be joined in pairs
$(x,y)$ in such a way that each $x$ is paired with exactly 
one $y$ and each $y$ with exactly one $x$.

\smallskip

b) The cardinality of $X$ is called ``less or equal''
to that $Y$, symbolically $|X|\le |Y|$, if there is a subset
$X^{\prime}\subset Y$ such that $|X|=|X^{\prime}|$. 

\smallskip

After these two {\it definitions}, the famous Cantor's {\it theorem}
can be proved in several lines:

\smallskip

c) The set of all subsets of $X$, symbolically $P(X)$,
has cardinality strictly larger than that of $X$.

\smallskip

Since we may iterate this construction, forming consecutively $P(P(X))$,
$P(P(P(X)))$, $\dots$, we see that 
{\it there exists an infinite scale of infinities of
growing sizes.}

\smallskip

The proof of c) consists of two remarks. The first
one says that $|X|\le |P(X)|$, because $X$ can be in a tautological
way paired with a part of $P(X)$ consisting of
one--element subsets of $X$. 

\smallskip

The second remark is
(a remake of) the famous Cantor's diagonal argument,
using {\it reductio ad absurdum}. Imagine that
$|X|=|P(X)|$. Then we can pair each $x\in X$ with some $S_x\subset X$
in such a way that {\it any} subset $S\subset X$ has the form
$S_y$ for some $y\in S$. Choose such a pairing (technically,
one--to--one correspondence). Define 

\smallskip

$S=$ {\it the set of all $x$ such that $x\notin S_x$}.

\smallskip
This $S$ must be of the form $S_y$ for some $y\in X$, but
then both logical possibilities, $y\notin S_y$ and $y\in S_y$
lead to a contradiction, so that the postulated
one--to--one correspondence cannot exist.

\smallskip

Of course, the last key argument goes back to the ancient
``liar's paradox''. It was revived again in a different context 
in the 20th century by Tarski and G\"odel. Tarski's theorem
features the ominous, at least for the purposes of this conference, 
name {\it ``inexpressibility of truth''}. 

\smallskip

In the final analysis, 
self--referentiality was used to produce several deep  mathematical arguments,
and this became possible only when the mathematical universe
became so extended that the language of mathematics could
be embedded into this universe as a part of it.
In particular, Leibniz's dream of merging language with meta--language
became a reality.

%\newpage

\bigskip

\centerline{\bf 3}

\medskip

\hfill{\it  The best test of truth is the power of the thought }

\hfill{\it to get itself accepted in the competition of the market.}

\smallskip

\hfill{\it Justice Oliver Wendell Holmes, Jr (1919)}

\medskip

When Cantor first presented his diagonal argument in a letter to
Dedekind in 1873, it was worded differently and used only
to prove that the cardinality of the natural numbers is
strictly less than that of the real numbers. The discovery of  the
proof itself was in a sense hardly more important than
the discovery of the definition of what it means, for one infinity to be
larger than another one.

\smallskip

As soon as this was achieved, Cantor started thinking about the
cardinality of the reals  compared with that of the
pairs of reals, or, geometrically, sets of points of
a curve and of a surface respectively. They turned out to be
equal! If we have a pair of numbers $(\alpha ,\beta )$ in $(0,1)$,
Cantor suggested to produce from them the third number $\gamma \in (0,1)$
by putting decimal digits of $\alpha$ to the odd places,
and those of $\beta$ to the even places. One sees, that
vice versa,  $(\alpha ,\beta )$ can be reconstructed from $\gamma$.
Dedekind, who was informed by Cantor's letter 
about this discovery as well, remarked that
this does not quite work because some rational numbers have
two decimal representations, such as $0,499999\dots =0,5000000\dots$.
Cantor had to spend some time to amend the proof, but
this was a minor embarrassment, in comparison with the
fascinating novelty of the fact itself:
{\it ``Ce que je vous ai communiqu\'e tout r\'ecemment est pour moi
si inattendue, si nouveau, que je ne pourrai pour ainsi dire pas arriver
\`a une certaine tranquillit\'e d'esprit avant que je n'aie re\c{c}u,
tr\`es honor\' e ami, votre jugement
sur son exactitude. Tant que vous ne m'aurez pas approuv\'e, je ne puis que dire: 
je le vois, mais je ne le crois pas''}, as Cantor famously wrote
to Dedekind.

\smallskip

This returns us to the basic question on the nature of truth.

\smallskip

We are reminded that the notion of ``truth'' is a reification
of a certain relationship between humans and {\it texts/utterances/statements}, the relationship
that is called ``belief'', ``conviction'' or ``faith'', and which
itself should be analyzed, together with other primary notions invoked in
this definition.

\smallskip

Professor Blackburn in [Bl] extensively discussed other relationships of humans
to texts, such as {\it scepticism, conservatism, relativism, deflationism.}
However, in the long range all of them are secondary in the practice of a researcher
in mathematics.

\smallskip

So I will return to truth.

\smallskip

I will skip analysis of the notion of ``humans'' :=) and will
only sketch what must be said about  texts, sources
of conviction, and methods of conviction peculiar to mathematics.

\smallskip

{\it Texts.} Alfred North Whitehead allegedly said that all of 
Western philosophy was but a footnote to Plato.

\smallskip

The underlying metaphor of such a statement is: ``Philosophy is a text'',
the sum total of all philosophic utterances.

\smallskip

Mathematics decidedly is {\it not} a text, at least not in the same
sense as philosophy. There are no authoritative
books or articles  to which subsequent 
generations turn again and again
for wisdom. Except for historians, nobody reads Euclid,
Newton, Leibniz or Hilbert in order to study
geometry, calculus or mathematical logic. 
The life span of any mathematical paper or book
can be years, in the best (and exceptional) case decades.
Mathematical wisdom, if not forgotten, lives as an invariant of all
its (re)presentations in a permanently self--renewing discourse. 

\smallskip
 
{\it Sources and methods  of conviction.} Mathematical truth is not revealed, 
and its acceptance is not imposed by any authority.

\smallskip

Moreover,  mathematical truth decidedly is not something that can be
ascertained, as Justice Oliver Wendell Holmes put it, by
``the majority vote of the nation that could lick all the others'' 
(quoted from [Pe]), or  by   acceptance ``in the competition of the market''.
In short, it is not a democratic value.

\smallskip

Ideally, the truth of a mathematical statement is ensured by {\it a proof},
and the ideal picture of a proof is a sequence of elementary arguments
whose rules of formation are explicitly laid down before the proof even begins,
and ideally are common for all proofs that have been devised
and can be devised in future. The admissible starting points
of proofs, ``axioms'', and terms in which they are formulated,
should also be discussed and made explicit.

\smallskip

This ideal picture is so rigid that it can itself become
the subject of mathematical study, which was actually performed 
and led to several
remarkable discoveries, technically all related to the effects of
merging language with metalanguage and self--referentiality.

\smallskip

Of course, the real life proofs are rendered in a peculiar mixture
of a natural language, formulas, motivations, examples.
They are much more condensed than imaginary formal proofs.
The ways of condensing them are not systematic in any way. We are prone to
mistakes, to taking on trust others' results that can be mistaken as well, and
to relying upon authority and revelations from our teachers.
(All of this should have been discussed together with the
notion of ``humans'' which I have wisely avoided.)

\smallskip

Moreover, the discovery of truth may, and usually does, involve
experimentation, nowadays vast and computer--assisted, false steps, sudden insights and all that which makes mathematical creativity
so fascinating for its adepts.

\smallskip

One metaphor of proof is a route, which might be
a desert track  boring and
unimpressive  until one finally reaches the oasis of one's destination, or
 a foot path in green hills, exciting and energizing, opening great vistas
of unexplored lands and seductive offshoots, leading far away even after 
the initial destination point has been reached.

\bigskip

\centerline{\bf 4}

\medskip

\hfill{\it [...] ``mismanagement and grief'':  here you have that}

\hfill{\it enormous distance between cause and effect covered in one line.}

\hfill{\it Just as math preaches how to do it.}

\smallskip

\hfill{\it J.~Brodsky. On ``September 1, 1939'' by W.~H.~Auden.}

\medskip

Mathematics is most visible to the general public  when it is posits itself
as an applied science, and in this role the notion of
mathematical truth acquires distinctly new features.
For example, our initial discussion of $\pi$
as an essentialy non--finitary (``irrational'') real number becomes pointless;
whenever $\pi$ enters any practical calculation, the first few digits are all that matters.  

\smallskip

In a wider context than just applied science, mathematics
can be fruitfully conceived as a toolkit containing powerful cognitive devices.
I have argued elsewhere ([Ma1], [Ma2]) that these devices can
be roughly divided into three overlapping domains: models,
theories,  and metaphors.  Quoting from [Ma2],

\bigskip

%\hangindent 15pt  

``A mathematical {\it model} describes a certain range of phenomena qualitatively
or quantitatively but feels uneasy pretending to be something more.

\smallskip

From Ptolemy's epicycles (describing planetary motions, ca 150) 
to the Standard Model (describing interactions of elementary particles, ca 1960), quantitative models cling to the observable reality
by adjusting numerical values of sometimes dozens of 
free parameters ($\ge 20$ for the Standard Model). Such models 
can be remarkably precise.

\smallskip

Qualitative models offer insights into {\it stability/instability},
{\it attractors} which are limiting states tending to occur independently
of initial conditions, {\it critical phenomena} in complex
systems which happen when the system crosses a boundary
between two phase states, or two basins
of different attractors.  [...]

\smallskip

What distinguishes a (mathematically formulated physical) {\it theory} 
from a model is primarily its 
higher aspirations. A modern physical
theory generally purports that it would describe the world
with absolute precision if only it (the world)
consisted of some restricted variety of
stuff: massive point particles obeying only the
law of gravity; electromagnetic field in a vacuum;
and the like.  [...]

\smallskip

A recurrent driving force generating theories is a concept
of a reality beyond and above the material world,
reality which may be grasped only by mathematical tools.
From Plato's solids to Galileo's ``language of nature''
to quantum superstrings, this psychological attitude
can be traced sometimes even if it conflicts with the
explicit philosophical positions of the researchers.   

\smallskip

 A (mathematical) {\it metaphor}, when it aspires
to be a cognitive tool, postulates that some 
complex range of phenomena might be compared to
a mathematical construction. The most recent
mathematical metaphor I have in mind is
Artificial Intelligence (AI). On the one hand, AI
is a body of knowledge related to computers and
a new, technologically created reality, consisting of
hardware, software, Internet etc. On the other hand,
it is a potential model of functioning of biological
brains and minds. In its entirety, it has not reached
the status of a model:
we have no systematic, coherent and extensive
list of correspondences between chips and neurons,
computer algorithms and brain algorithms.
But we can and do use our extensive knowledge of 
algorithms and computers
(because they were created by us) to generate
educated guesses about structure and function of the central
neural system [...].

\smallskip

A mathematical theory is an invitation 
to build applicable models. A mathematical metaphor
is an invitation to ponder upon what we know.''

\bigskip

As an aside, let us note that George Lakoff's {\it definition}
of poetic metaphors such as ``love is a journey'' 
in [La]  is itself expressed as a mathematical metaphor
using the characteristic Cantor--Bourbaki mental images and vocabulary:
{\it ``More technically, the metaphor can be understood as a mapping
(in the mathematical sense) from a source domain (in this case, journeys)  to a target domain
(in this case, love). The mapping is tightly structured. There are ontological
correspondences, according to which entities in the domain of love
(e.~g. the lovers, their common goals,
their difficulties, the love relationship, etc.) correspond systematically
to entities in the domain of a journey (the travellers, the vehicle,
destinations, etc.).}''

\smallskip

When a mathematical construction  is used as a cognitive tool,  the 
discussion of truth becomes loaded with  new meanings:
a model, a theory or a metaphor must be
true to a certain reality, more tangible and real than the Platonic
``reality''  of pure mathematics.  In fact, philosophers of science routinely discussed
truth precisely in this context. Karl Popper's vision of
scientific theories in terms  of falsifiability (vs  verifiability)
is quite appropriate in the context of highly mathematicised theories as well.

\smallskip

What I want to stress here, however, is one aspect of  contemporary 
mathematical models which is historically very recent. 
Namely,  models are more and more widely used as ``black boxes''
with hidden computerized input
procedures, and oracular outputs prescribing behavior of 
human users.  

\smallskip

Mary Poovey, discussing from this viewpoint financial markets, remarks
in her insightful essay [Po] that what she calls ``representations'',
basically computerized bookkeeping or the numbers a trader enters in a computer,
tend to replace the actual exchange of cash or commodities.
``This conflation of  representation and exchange has all kinds of material effects,
[...] for when representation can influence or take the place of exchanges,
the values at stake become notional too: they can grow exponentially or collapse
at the stroke of key''.

\smallskip

In fact, actions of traders, banks, hedge funds and alike are  to a considerable degree determined
by the statistical models of financial markets encoded in the software of their computers.
These models thus become a hidden and highly influential
part of the actions,  our  computerized  ``collective unconscious''.
As such, they cannot even be judged according to the usual criteria
of choosing models which better reflect the behavior of a process being modeled.
They are part of any such process.

\smallskip

What becomes more essential than their empirical adequacy,  is, for example, their stabilizing
or destabilizing potential.  Risk management assuming mild variability and
small risks can collapse when a disaster occurs, ruining many participants of
the game;  risk management  based upon models that use pessimistic ``L\'evy distributions''
rather than omnipresent Gaussians paradoxically tends to flatten the shock waves
and thus to avoid major disasters (cf.  [MandHu]).

\newpage

\centerline{\bf 5}

\medskip

\hfill{\it There have been dramatic changes in the way in which}

\hfill{\it the motion of the crowd is modeled in recent years.}

\smallskip

\hfill{\it R.~Clemens, R.~Hughes, in [ClHu].}

\medskip

When in the 20th century mathematicians
got involved in heated discussions
about the so called ``Crisis in Foundations of Mathematics'',
several issues were intermingled.

\smallskip

Philosophically--minded logicians and professional philosophers
were engaged with the nature and accessibility of mathematical
truth (and reliability of
our mental tools used in the process of acquiring it).

\smallskip

Logicists (finitists, formalists, intuitionists) were elaborating severe
normative prescriptions trying to outlaw dangerous mental
experiments with infinity, non--constructivity and {\it reductio ad absurdum.}

\smallskip

For a working mathematician, when he/she is concerned 
at all, ``foundations''  is simply a general term for the historically variable
set of rules and principles of organization of the body of
mathematical knowledge, both existing and being created.
From this viewpoint, the most influential foundational achievement
in the 20th century was an ambitious project
of the Bourbaki group, building all mathematics, including
logic, around set--theoretical ``structures'' and making Cantor's
language of  sets a common 
vernacular of algebraists, geometers, probabilists and all other
practitioners of our trade. These days, this vernacular,
with all its vocabulary and ingrained mental habits,
is being slowly replaced by the languages of category theory
and homotopy theory and their higher extensions. Respectively,
the basic ``left--brain'' intuition of sets, composed of distinguishable
elements, is giving way to a new, more
``right brain'' basic intuition dealing with space--like and
continuous
primary images, both deformable and deforming.

\smallskip

In the Western ethnomathematics, truth is best understood 
as a central value, ever to be pursued, rather than anything
achieved.  Practical efficiency, authority, success in competition,
faith, all other clashing values must recede in the mind
of a mathematician when he or she sets down to do their job. 

\smallskip

The most interesting intracultural interactions of mathematics
such as symbolized by this conference 
are as well those that are not direct but
rather proceed with the mediation of value systems.

\newpage

\centerline{\bf Coda}

\medskip

Every four years, mathematicians from all over the world meet
at the International Congresses (ICM), to discuss whatever
interesting developments happened recently in their
domains of expertise. One of the traditions
of these Congresses is a series of lectures for general public.

\smallskip

In 1998, our Congress met in Berlin, and Hans Magnus Enzensberger,
the renowned poet and essayist, deeply interested in
mathematics, spoke about ``Zug\-br\"ucke au{\ss}er Betrieb: die Mathematik
im Jenseits der Kultur'': the
drawbridge to the castle of mathematics is out of service. 
The main concern of his talk was a deplorable lack of mathematical
culture and communication between the general public and mathematicians,
leading to alienation and mutual mistrust.

\smallskip

At the end of his talk ([Enz]) Enzensberger quotes an imaginary
dialogue from [St], where a mathematician is chatting 
with a fictitional layman ``Seamus Android''.

\medskip

``Mathematician: It's one of the most important discoveries
of the last decade!

\smallskip

Android: Can you {\it explain} it in words ordinary mortals
can understand?

\smallskip

Mathematician: Look, buster, if ordinary mortals
could understand it, you wouldn't need mathematicians
to do the job for you, right? You can't get a feeling for what's going on
without understanding the technical details. How can I talk about manifolds
without mentioning that the theorems only work if the
manifolds are finite--dimensional para--compact Hausdorff
with empty boundary?

\smallskip

Android: Lie a bit.

\smallskip

Mathematician: Oh, but I couldn't do that!

\smallskip

Android: Why not? Everybody {\it else} does.''

\medskip

And here I must play God and say to both Android
and Mathematician: ``Oh, no! Don't lie --- because
everybody else does.''
  
\bigskip

\centerline{\bf References}

\medskip

[Bl] S.~Blackburn.  {\it  Truth and Ourselves: the Elusive  Last Word.}
Keynote talk at the Balzan Symposium ``Truth'', May 2008.

\smallskip

[BoHa] L.~Bovens, S.~Hartmann. {\it Bayesian Epistemology.}
Clarendon Press, Oxford, 2003.

\smallskip

[ClHu]  R.~C.~Clemens, R.~L.~Hughes. {\it Mathematical Modelling of
a Mediaeval Battle: the Battle of Agincourt, 1415.} Math. and Computers
in Simulation, 64:2 (2004), 259--269.

\smallskip

[Dau] J.~W.~Dauben. {\it Georg Cantor. His Mathematics and Philosophy of the Infinite.}  Princeton University Press, 1990.

\smallskip

[DavHe] P.~Davis, R.~Hersh. {\it The Mathematical
Experience.} Birkh\"auser Boston, 1986.

\smallskip

[Enz] H.~M.~Enzensberger. {\it Drawbridge Up. Mathematics -- a Cultural
Anathema.} A.~K.~Peters, Natick, Mass., 1999.

\smallskip

[La] G.~Lakoff. {\it  The Contemporary Theory of Metaphor.} In: A. Ortony (ed.),
Metaphor and Thought (2nd ed.). Cambridge Univ. Press, 1993.

\smallskip

[MandHu]  B.~Mandelbrot, R.~Hudson. {\it The (Mis)behavior of Markets: a Fractal View of
Risk, Ruin and Rewards.} Profile, 2005.

\smallskip

[Ma1] Yu.~Manin.  {\it Mathematics as Metaphor.} (Selected Essays, with Foreword by
F.~Dyson). American Math. Society, 2007.

\smallskip

[Ma2]  Yu.~Manin. {\it Mathematical Knowledge: Internal, Social and Cultural Aspects. }
(Introductory Chapter to vol.~2 of  ``La Matematica'', Einaudi,
ed. by C.~Bartocci and P.~Odifreddi, reproduced in [Ma1], pp. 3--26).  Preprint math.HO/0703427

\smallskip

[Pe] J.~D.~Peters. {\it Courting the Abyss: Free Speech and the Liberal Tradition.} Chicago, 2005.

\smallskip

[Po] M.~Poovey. {\it Can Numbers Ensure Honesty?
Unrealistic expectations and the US accounting scandal.}
Notices of the AMS, vol. 50:1, Jan. 2003, pp. 27--35.

\smallskip

[St] I.~Stewart. {\it The Problems of Mathematics.} Oxford Univ. Press, 1987.

\enddocument